\documentclass[11pt,twoside]{amsart}
\usepackage{amsmath, amsthm, amscd, amsfonts, amssymb, graphicx, color}
\usepackage{blkarray}
\usepackage{caption}
\usepackage{subcaption}
\usepackage[colorlinks=true,linkcolor=black,anchorcolor=black,citecolor=black,filecolor=black,menucolor=black,runcolor=black,urlcolor=black]{hyperref}
\newtheorem{thm}{Theorem}[section]

\newtheorem{lem}[thm]{Lemma}

\newtheorem{defn}[thm]{Definition}

\numberwithin{equation}{section}
\newtheorem{ex}{Example}
\numberwithin{equation}{thm}
\def\pn{\par\noindent}

\begin{document}

\title{On the Spectral properties of Andr\'asfai Graphs}
\author{Bharani Dharan K and Radha S$^*$}

\thanks{{\scriptsize
\hskip -0.4 true cm MSC(2010): Primary: 05C12, 05C25, 05C50
Received: dd mmmm yyyy, Accepted: dd mmmm yyyy.\\
$*$Corresponding author\: Radha S\\
Mail.Id.: radha.s@vit.ac.in}}
\maketitle

\begin{center}
Communicated by\; Radha S
\end{center}

\begin{abstract} 
In this paper, we investigate the spectral properties of Andr\'asfai graphs, focusing on key parameters: the second-largest and smallest eigenvalues, the number of distinct eigenvalues, and the multiplicities of the eigenvalues 1 and -1. The results obtained reveal insights into the connectivity, the structural properties, and the spectral distinctiveness. 
\end{abstract}

\vskip 0.1 true cm



\section{\bf Introduction}
\vskip 0.1 true cm
Spectra of some graphs related to networks are intensively useful for identifying drugs for complex diseases \cite{19}, investigating global commerce networks \cite{17}, determining the stability of a system \cite{20}, and many other processes.  Most recent studies have concentrated on the relationships between graphs and the count of distinct eigenvalues of the matrices associated with them. The concept of the distinct eigenvalues of the matrices associated with the graph is studied in works like \cite{4,5,6,7,8,9,14}. The count of the distinct eigenvalues of the perturbation matrices and a bound for the same are discussed in the recent works \cite{10,11,12,13}. \\
The solution to the K$\ddot{o}$nigsberghg bridge problem \cite{15,16} and the family of Cayley graphs has found amazing applications in many real-world systems. The theory of networks\cite{27}, which investigates complicated interacting units represented as graphs, makes extensive use of the family of Cayley graphs. It has been shown that information about the network's structural properties as well as the dynamic behavior of the connected complex system may be found in the adjacency matrix of the network's spectra. The largest eigenvalue encapsulates significant information and correlates with the entrainment of the diffusely coupled dynamical units in the network \cite{21}. Additionally, sparse symmetric matrices \cite{22} have provided a technique for calculating the statistical properties of the second largest eigenvalue and the constituents of the associated characteristic vector. It should be noted that many network structural features are computationally difficult to establish; nevertheless, spectral measures often provide valuable insights into the structure of networks and are computationally more straightforward to compute. For instance, computing several network growth aspects is computationally difficult. Fortunately, the second-largest eigenvalue, which can be calculated fast (in $O(n^3)$, where n is the number of network nodes), is strongly connected to these characteristics\cite{23}.\\
The majority of research on the second largest eigenvalue is done for some arbitrary regular graphs, with infrequent studies conducted for other networks. All of these investigations opine that the second largest eigenvalue of a graph can provide useful information into the characteristics of the underlying network topology\cite{25}. In particular, the second largest eigenvalue of a graph determines whether a network is appropriate for a certain application, and the minimal value of the second largest eigenvalue is frequently desired \cite{26}.\\
In the second section of this paper, we provide some basic definitions and findings relevant to our work. In section 3, we have stated and proved our main results on the spectrum of Andr\'asfai graphs $And(k)$, such as the number of distinct eigenvalues, the least eigenvalue, and the second largest eigenvalue.
 \section{Preliminaries}
	       Let G be a graph on n vertices and let $A_G$ be the adjacency matrix of the graph G, then the eigenvalues of $A_G$ is usually denoted by $\lambda_0,\lambda_{1},\lambda_{2},\lambda_{3},.......,\lambda_{n-1}$ where $\lambda_0\geq\lambda_{1}\geq\lambda_{2}\geq\lambda_{3}\geq........\geq\lambda_{n-1}$ (i.e.) $\lambda_0,~\lambda_1,~and~\lambda_{n-1}$ are the largest, second largest and smallest eigenvalues of $A_G$ respectively \cite{3}.  In our paper, we denote the eigenvalues of $A_{And(k)}$ by $x_0,x_1,x_2,...,x_{n-1}$ where $x_i's$ need not be of the form $x_0\geq x_1\geq x_2 \geq \dots \geq x_{n-1}$ (i.e.) $x_i's$ need not be equal to $\lambda_i$ where $0\leq i\leq n-1$.
	\begin{defn}
		Let $k\geq1$ be any natural number and take $n=3k-1$.  Andr\'asfai Graph is a Cayley graph over the additive group $\mathbb{Z}_n$,(i.e.) Cay$(\mathbb{Z}_{n},S)$ where the generating set S=$\{$ x $|$  x$\in\mathbb{Z}_{n}$ and x$\equiv$1 mod 3 $\}$ is a subset of $\mathbb{Z}_{n}$. Generally, Andr$\acute{a}$sfai graphs are denoted by $\it{And}(k)$. \\For example, $And(3)$, $And(4)$ and $And(5)$ are shown in the figures \ref{And(3)},\ref{And(4)} and \ref{And(5)} given below. 
\begin{figure}[htp]
    \begin{subfigure}{.3\textwidth}
        \includegraphics[width=1\linewidth]{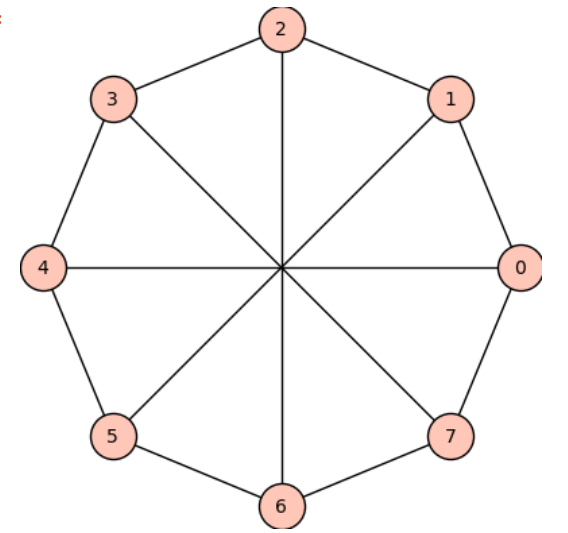}  
        \caption{And(3)}
        \label{And(3)}
    \end{subfigure}
    \begin{subfigure}{.3\textwidth}
        \centering
        \includegraphics[width=1\linewidth]{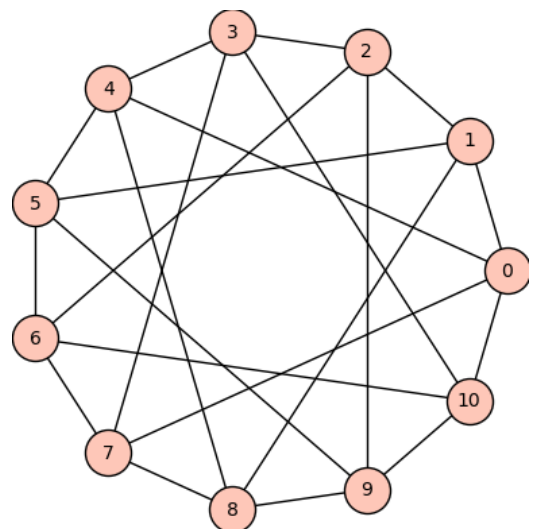}  
        \caption{And(4)}
        \label{And(4)}
    \end{subfigure}
    \begin{subfigure}{.3\textwidth}
        \centering
        \includegraphics[width=1\linewidth]{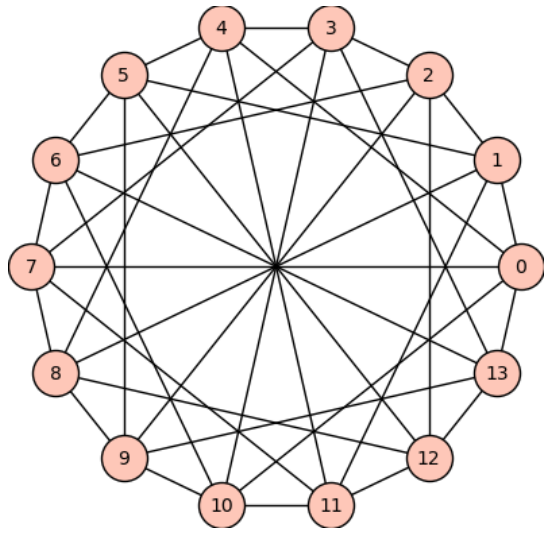}  
        \caption{And(5)}
        \label{And(5)}
    \end{subfigure}
    \caption{Andr\'asfai Graphs for $k=3,4,5$}
    \label{FIGURE LABEL}
\end{figure}

\end{defn}
	Let  G=$\it{And}(k)$, k$\geq$2. Then G must be k-regular, circulant and  the adjacency matrix of G is
 \[
A_{And(k)}=
\begin{blockarray}{ cccccccccccccc}
 &0&1&2&3&4&5&.&.&.&3k-5&3k-4&3k-3&3k-2 \\
\begin{block}{c(ccccccccccccc)}
0&0&1&0&0&1&0&.&.&.&1&0&0&1\\
1&1&0&1&0&0&1&.&.&.&0&1&0&0\\
2&0&1&0&1&0&0&.&.&.&0&0&1&0\\
3&0&0&1&0&1&0&.&.&.&1&0&0&1\\
4&1&0&0&1&0&1&.&.&.&0&1&0&0\\
5&0&1&0&0&1&0&.&.&.&0&0&1&0\\
6&0&0&1&0&0&1&.&.&.&1&0&0&1\\
7&1&0&0&1&0&0&.&.&.&0&1&0&0\\
:&:&:&:&:&:&:&:&:&:&:&:&:&:\\
:&:&:&:&:&:&:&:&:&:&:&:&:&:\\
:&:&:&:&:&:&:&:&:&:&:&:&:&:\\
3k-5&1&0&0&1&0&0&.&.&.&0&1&0&0\\
3k-4&0&1&0&0&1&0&.&.&.&1&0&1&0\\
3k-3&0&0&1&0&0&1&.&.&.&0&1&0&1\\
3k-2&1&0&0&1&0&0&.&.&.&0&0&1&0\\
\end{block}
\end{blockarray}
\]
Then the eigenvalues $x_{l}$'s of the above matrix $A_{And(k)}$ are given by
	\begin{equation}\label{equ1}
		x_{l}=\sum_{j=0}^{3k-2}a_{j}\omega^{lj},\hspace{1cm}0\leq l \leq 3k-2
	\end{equation}
	where $a_{j}$ is the $j^{th}$ entry of the first row of $A_{And(k)}$ and $\omega$ is the $(3k-1)^{th}$ root  of unity.
	\\From, \eqref{equ1}  we have
	\begin{equation}\label{equ2}
		x_{l}=\omega^{l}+\omega^{4l}+.....+\omega^{3k-5}+\omega^{3k-2}\\
	\end{equation}
	\begin{equation}   x_{l}=\omega^{l}+\omega^{4l}+.....+\omega^{-4l}+\omega^{-l}\nonumber 
	\end{equation}
	Case(i): When k is even
 \begin{equation*}
	\begin{aligned}
            x_{l}=&\sum_{j=0}^{\frac{k-2}{2}}\left(\omega^{(3j+1)l}+\omega^{-(3j+1)l}\right)\\ \nonumber
            x_{l}=&\sum_{j=0}^{\frac{k-2}{2}}\left(e^{\frac{2(3j+1)l\pi i}{n}}+e^{\frac{2(3j+1)l\pi i}{n}}\right)\\ \nonumber
	\end{aligned}
 \end{equation*}
	\begin{equation}\label{eigenwhenkiseven}
		x_{l}=2 \left[\sum_{j=0}^{\frac{k-2}{2}} cos\left(\frac{2(3j+1)l\pi}{n}\right)\right]
	\end{equation}\\ 
	Case(ii): When k is odd
 \begin{equation*}
	\begin{aligned}
            x_{l}=&\sum_{j=0}^{\frac{k-3}{2}}\left(\omega^{(3j+1)l}+\omega^{-(3j+1)l}\right)+(-1)^l\\ \nonumber
            x_{l}=&\sum_{j=0}^{\frac{k-3}{2}}\left(e^{\frac{2(3j+1)l\pi i}{n}}+e^{\frac{2(3j+1)l\pi i}{n}}\right)+(-1)^l\\ \nonumber
	\end{aligned}
 \end{equation*}
	\begin{equation}\label{eigenwhenkisodd}
		x_{l}=2 \left[\sum_{j=0}^{\frac{k-3}{2}} cos\left(\frac{2(3j+1)l\pi}{n}\right)\right]+(-1)^{l}
	\end{equation}
	\eqref{eigenwhenkiseven} and \eqref{eigenwhenkisodd} are the two structures of the eigenvalues of $A_{And(k)}$.\\
 \section{Main Results}\label{sec2}
	\begin{thm} \label{theorem1}
		Let G be an Andr\'asfai Graph $And(k)$.  Then G has $k+\lceil{\frac{k}{2}}\rceil$ distinct adjacency eigenvalues.
	\end{thm}
	\begin{proof}
		Before proving the main result, first, we will prove that,
		\begin{center}$x_{i}=x_{3k-1-i}$
		\end{center}
		We know that,
		\begin{equation*}
			\begin{aligned}
				x_{i} & =\omega^{i}+\omega^{4i}+\omega^{7i}+......+\omega^{(3k-5)i}+\omega^{(3k-2)i}\\
				& =\omega^{i}+\omega^{4i}+\omega^{7i}+......+\omega^{-7i}+\omega^{-4i}+\omega^{-i}\\
				& =\omega^{3k-1+i}+\omega^{3k-1+4i}+\omega^{3k-1+7i}+......+\omega^{3k-1-7i}+\omega^{3k-1-4i}+\omega^{3k-1-i}\\
				& =x_{3k-1-i}.\\
			\end{aligned}
		\end{equation*}
		\\Hence,\begin{equation}
			x_{i}=x_{3k-1-i}\label{(1)}
		\end{equation}
		\\Case(i): k is even\\
		When k is even, n=3k-1 is odd. Using \eqref{(1)}, $A_{And(k)}$ has an eigenvalue $x_{0}$ with multiplicity equal to 1.\\
		To find the multiplicities of other eigenvalues:\\
		Suppose that \\
		\begin{equation*}
			\begin{aligned}
				 x_{l}=x_{m},\text{ for }l\neq m \text{ and } 1\leq l,m\leq3k-2
			\end{aligned}
		\end{equation*} 
		\\Then,
		\begin{equation*}
			\begin{aligned}
				x_{l}-x_{m}=0
			\end{aligned}
		\end{equation*}
  \begin{equation*}
      \implies 2\left[\sum_{j=0}^{\frac{k-2}{2}} cos\left(\frac{2(3j+1)l\pi}{n}\right)\right]-2\left[\sum_{j=0}^{\frac{k-2}{2}} cos\left(\frac{2(3j+1)m\pi}{n}\right)\right]=0
  \end{equation*}
		\begin{equation*}
			\implies 2 \left[\sum_{j=0}^{\frac{k-2}{2}} cos\left(\frac{2(3j+1)l\pi}{n}\right)-\sum_{j=0}^{\frac{k-2}{2}}cos\left(\frac{2(3j+1)m\pi}{n}\right)\right]=0
		\end{equation*}
            \begin{equation}\label{main1}
			\implies 4 \left[\sum_{j=0}^{\frac{k-2}{2}}\left(sin\left(\frac{2(3j+1)(l+m)\pi }{n}\right)sin\left(\frac{2(3j+1)(m-l)\pi}{n}\right)\right) \right]=0
		\end{equation}
	
		For  $l\neq m$  the above equation \eqref{main1} is true only when $m=-l$\\ where $m$ and $l$ are the elements of the group $\mathbb{Z}_{n}$.
		Then it can be written as $l=-m$ mod n,\\
		(i.e) $m=3k-1-l$\\
		Hence, $x_{l}=x_{m}$ only when $m=3k-1-l$.\\
		Then all other eigenvalues have multiplicity equal to 2, except $x_{0}$.\\
		\\Case(ii): When k is odd\\
		For an odd k, n is even.  Using \ref{(1)}, $A_{And(k)}$ has two eigenvalues $x_{0}$ and $x_{\frac{3k-1}{2}}$ with multiplicities equal to 1. To find the multiplicities of other eigenvalues:
		\\Suppose that 
		\begin{equation*}
			\begin{aligned}
				 x_{l}=x_{m},\text{ for }l\neq m
			\end{aligned}
		\end{equation*}
		Then,
		\begin{equation*}
			\begin{aligned}
				x_{l}-x_{m}=0
			\end{aligned}
		\end{equation*}
		\begin{multline*}
			2 \left[\sum_{j=0}^{\frac{k-3}{2}} cos\left(\frac{2(3j+1)l\pi}{n}\right)\right]+(-1)^{l}-2 \left[\sum_{j=0}^{\frac{k-3}{2}}cos\left(\frac{2(3j+1)m\pi}{n}\right)\right]-(-1)^{m}=0
		\end{multline*}
             \begin{multline*}
			2 \left[\sum_{j=0}^{\frac{k-3}{2}} \left(cos\left(\frac{2(3j+1)l\pi}{n}\right)-cos\left(\frac{2(3j+1)m\pi}{n}\right)\right)\right]+(-1)^{l}-(-1)^{m}=0
		\end{multline*}
            \begin{equation}\label{main2}
			\implies 4 \left[\sum_{j=0}^{\frac{k-3}{2}}\left(sin\left(\frac{2(3j+1)(l+m)\pi }{n}\right)sin\left(\frac{2(3j+1)(m-l)\pi}{n}\right)\right) \right]]+(-1)^{l}-(-1)^{m}=0
		\end{equation}
		
		For $l\neq m$ the above equation \eqref{main2} becomes true only when $m=-l$\\ where $l$ and $m$ are the elements of the group $\mathbb{Z}_{n}$.
		Then it can be written as $l=-m$ mod n,\\
		(i.e) $m=3k-1-l$\\
		Hence, $x_{l}=x_{m}$ only when $m=3k-1-l$\\
		Then all other eigenvalues have multiplicity equal to $2$ except $x_{0}$ and $x_{\frac{3k-1}{2}}$.\\
		\\From the above discussion, \\ \hspace{1cm}when k is even, the number of distinct eigenvalues of $A_{G}$ is equal to $\frac{3k+1}{2}$,\\ \hspace{1cm} when k is odd, the number of distinct eigenvalues of $A_{G}$ is equal to $\frac{3k}{2}$.\\
		\\Without loss of generality we can say the number of different eigenvalues of $A_{G}$ is equal to $k+\lceil\frac{k}{2}\rceil$.
	\end{proof}
	\begin{thm} \label{theorem2}
		For the adjacency matrix of the Andr\'asfai graph $And(k)$, $x_{k}=x_{2k-1}$ is the smallest eigenvalue.
	\end{thm}
	\begin{proof}
		The theorem can be proved by the method of contradiction.\\
		Suppose that there exists $l\neq 0,k,2k-1$, such that $x_{l}<x_{k}$, then
		\begin{equation}\label{leastequation}
			x_{l}-x_{k}<0
		\end{equation}
		Case(i): When k is even, \eqref{leastequation} becomes
		\begin{multline*}
			\implies 2 \left[\sum_{j=0}^{\frac{k-2}{2}} cos\left(\frac{2(3j+1)l\pi}{n}\right)\right]-2 \left[\sum_{j=0}^{\frac{k-2}{2}}cos\left(\frac{2(3j+1)k\pi}{n}\right)\right]<0
		\end{multline*}
 
		\begin{equation*}
			\implies 2 \left[\sum_{j=0}^{\frac{k-2}{2}} cos\left(\frac{2(3j+1)l\pi}{n}\right)-\sum_{j=0}^{\frac{k-2}{2}}cos\left(\frac{2(3j+1)k\pi}{n}\right)\right]<0
		\end{equation*}
            \begin{equation}\label{leastmain1}
			\implies 4 \left[\sum_{j=0}^{\frac{k-2}{2}}\left(sin\left(\frac{2(3j+1)(l+k)\pi }{n}\right)sin\left(\frac{2(3j+1)(k-l)\pi}{n}\right)\right) \right]<0
		\end{equation}
		
		The L.H.S on the inequality\eqref{leastmain1} is always non-negative for every $l\in\mathbb{Z}_{3k-1}$ and contradicts \eqref{leastmain1}.
		\\Hence, $x_{k}$ is the smallest eigenvalue of the adjacency matrix of a Andr\'asfai graph $And(k)$ when k is even.\\
		Case(ii): When k is odd, \eqref{leastequation} becomes 
  \begin{multline*}
			\implies 2 \left[\sum_{j=0}^{\frac{k-3}{2}} cos\left(\frac{2(3j+1)l\pi}{n}\right)\right]+(-1)^{l}-2 \left[\sum_{j=0}^{\frac{k-3}{2}}cos\left(\frac{2(3j+1)k\pi}{n}\right)\right]-(-1)^{k}<0
		\end{multline*}
             \begin{multline*}
			\implies 2 \left[\sum_{j=0}^{\frac{k-3}{2}} \left(cos\left(\frac{2(3j+1)l\pi}{n}\right)-cos\left(\frac{2(3j+1)k\pi}{n}\right)\right)\right]+(-1)^{l}-(-1)^{k}<0
		\end{multline*}
            \begin{equation}\label{leastmain2}
			\implies 4 \left[\sum_{j=0}^{\frac{k-3}{2}}\left(sin\left(\frac{2(3j+1)(l+k)\pi }{n}\right)sin\left(\frac{2(3j+1)(k-l)\pi}{n}\right)\right) \right]]+(-1)^{l}-(-1)^{k}<0
		\end{equation}
		
		The L.H.S of the inequalilty\eqref{leastmain2} is always non-negative for every $l\in \mathbb{Z}_{3k-1}$ and contradicts \eqref{leastmain2}.
		\\Hence, (by \eqref{(1)}) $x_{k}=x_{2k-1}$ is the smallest eigenvalue of the adjacency matrix of a Andr\'asfai graph $And(k)$.
		Hence proved.
		
	\end{proof}
	\begin{thm}
		For the adjacency matrix of the Andr\'asfai graph $And(k)$, $x_{k-1}=x_{2k}$ is the second largest eigenvalue.
		\begin{proof}
			The theorem can be proved by the method of contradiction.
			We know that $x_{0}$ is the greatest eigenvalue as the graph is k-regular. So, we consider the eigenvalues other than $x_{0}$.\\
			Suppose that there exists $l\neq0,k-1,2k$, such that, $x_{l}>x_{k-1}$,
			\begin{equation}\label{secondlargesteqa}
				x_{l}-x_{k-1}>0
			\end{equation}
			Case(i): When k is even, \eqref{secondlargesteqa} becomes 
   \begin{multline*}
			\implies 2 \left[\sum_{j=0}^{\frac{k-2}{2}} cos\left(\frac{2(3j+1)l\pi}{n}\right)\right]-2 \left[\sum_{j=0}^{\frac{k-2}{2}}cos\left(\frac{2(3j+1)(k-1)\pi}{n}\right)\right]>0
		\end{multline*}
 
		\begin{equation*}
			\implies 2 \left[\sum_{j=0}^{\frac{k-2}{2}} cos\left(\frac{2(3j+1)l\pi}{n}\right)-\sum_{j=0}^{\frac{k-2}{2}}cos\left(\frac{2(3j+1)(k-1)\pi}{n}\right)\right]>0
		\end{equation*}
            \begin{equation}\label{secondlargestmain1}
			\implies 4 \left[\sum_{j=0}^{\frac{k-2}{2}}\left(sin\left(\frac{2(3j+1)(l+k-1)\pi }{n}\right)sin\left(\frac{2(3j+1)(k-1-l)\pi}{n}\right)\right) \right]>0
		\end{equation}
			
			The L.H.S of the inequality \eqref{secondlargestmain1} is always non-positive for every  $l\in \mathbb{Z}_{3k-1}$ and contradicts \eqref{secondlargestmain1}.\\ Hence, $x_{k-1}=x_{2}$ is the second largest eigenvalue of the adjacency matrix of the Andr\'asfai graph $And(k)$ when k is even.
			\\Case(ii): When k is odd,\eqref{secondlargesteqa} becomes 
            \begin{multline*}
			\implies2 \left[\sum_{j=0}^{\frac{k-3}{2}} cos\left(\frac{2(3j+1)l\pi}{n}\right)\right]+(-1)^{l}-2 \left[\sum_{j=0}^{\frac{k-3}{2}}cos\left(\frac{2(3j+1)(k-1)\pi}{n}\right)\right]-(-1)^{k-1}>0
		\end{multline*}
             \begin{multline*}
			\implies2 \left[\sum_{j=0}^{\frac{k-3}{2}} \left(cos\left(\frac{2(3j+1)l\pi}{n}\right)-cos\left(\frac{2(3j+1)(k-1)\pi}{n}\right)\right)\right]+(-1)^{l}-(-1)^{k-1}>0
		\end{multline*}
            \begin{equation}\label{secondlargestmain2}
			\implies 4 \left[\sum_{j=0}^{\frac{k-3}{2}}\left(sin\left(\frac{2(3j+1)(l+k-1)\pi }{n}\right)sin\left(\frac{2(3j+1)(k-1-l)\pi}{n}\right)\right) \right]]+(-1)^{l}-(-1)^{k-1}>0
		\end{equation}
			
			The L.H.S of the above inequality is always non-positive for every $l\in\mathbb{Z}_{3k-1}$ and contradicts \eqref{secondlargestmain2}. \\ Hence, (By \eqref{(1)})$x_{k-1}=x_{2k}$ is the second largest eigenvalue of the adjacency matrix of the Andr\'asfai graph $And(k)$.
			Hence proved.
		\end{proof}
	\end{thm}
 
        \begin{lem}\label{-1eigenlemma}
            For an $And(k)$, -$1$ is an eigenvalue with multiplicity $1$ if and only if $k$ is odd.
        \end{lem}
        \begin{proof}
            Let us take $And(k)$ with $k$ as odd, which implies that 
            \begin{equation*}
             x_{l}=2 \left[\sum_{j=0}^{\frac{k-3}{2}} cos\left(\frac{2(3j+1)l\pi}{n}\right)\right]+(-1)^{l}
             \end{equation*}
             Substitute, $l=\frac{n}{2} \implies$
             \begin{equation*}
            \begin{aligned}
                x_{\frac{n}{2}}&=2 \left[\sum_{j=0}^{\frac{k-3}{2}} cos\left(\frac{2(3j+1)\left(\frac{n}{2}\right)\pi}{n}\right)\right]+(-1)^{\frac{n}{2}}\\
                x_{\frac{n}{2}}&=2 \left[\sum_{j=0}^{\frac{k-3}{2}} \cos\left((3j+1)\pi\right)\right]+(-1)^{\frac{n}{2}}
            \end{aligned}
            \end{equation*}
            When $k=1~mod~4$, $\frac{k-3}{2}$ is even and  $\frac{n}{2}$ is odd which makes $x_{\frac{n}{2}}$ is equal to minus one. 
            When $k=3~mod~4$,  $\frac{k-3}{2}$ is odd and $\frac{n}{2}$ is even which makes $x_{\frac{n}{2}}$ is equal to minus one.\\
            We can also prove the above necessary part by alternate method.
            When $k$ is odd, $n$ is even and there exists a integer $l=\frac{n}{2}$, such that the corresponding eigenvalue $x_l$ can be written as 
            \begin{equation*}
            \begin{aligned}
                x_l&=\omega^l+\omega^{4l}+\omega^{7l}+\dots+\omega^{(3k-2)l}\\
                x_{\frac{n}{2}}&=\omega^{\frac{n}{2}}+\omega^0+\omega^{\frac{n}{2}}+\dots+\omega^{\frac{n}{2}}+\omega^0+\omega^{\frac{n}{2}}\\
                &=\left(\frac{k+1}{2}\right)\omega^{\frac{n}{2}}+\left(\frac{k-1}{2}\right)\omega^0\\
                x_{\frac{n}{2}}&=-1
            \end{aligned}
            \end{equation*}
            
            Now, we need to prove the sufficient part (ie) if -1 is an eigenvalue of $And(k)$, then $k$ is odd. 
            Let us assume that -1 is an eigenvalue of $And(k)$, then
            \begin{equation*}
            \begin{aligned}
                x_l&=-1\\
                \sum_{j=0}^{k-1}\left(\omega\right)^{\left(3j+1\right)l}&=-1\\
                \sum_{j=0}^{k-1}\left(\omega\right)^{\left(3j+1\right)l}+\omega^0&=0
                \end{aligned}
            \end{equation*}
        
            \begin{equation}\label{minusoneeigen}              \left(\omega^l\right)^0+\left(\omega^l\right)^1+\left(\omega^l\right)^4+\left(\omega^l\right)^7+\dots+\left(\omega^l\right)^{3k-2}=0
            \end{equation}             
            There are $k+1$ terms in the equation\eqref{minusoneeigen} and each term is the $g^{th}$ root of unity, where $g=gcd(n,k+1)$.  If the  equation\eqref{minusoneeigen} becomes true then there exist $l\in \mathbb{Z}_{3k-1}$ such that, $l=\frac{n}{gcd(n,k+1)}$, where $gcd(n,k+1)>1$. Then
            \begin{equation}\label{gcd}
            gcd(n,k+1)=a,\text{ where }a>1
            \end{equation}
            Equation\eqref{gcd} implies that
            \begin{equation}\label{s1s2}
            \begin{aligned}
                3k-1&=as_1\\
                k+1&=as_2
            \end{aligned}   
            \end{equation}
            where $s_1$ and $s_2$ are positive integers and their values are dependent  on $k$ and $a$.
            Solving \eqref{s1s2}, we get
            \begin{equation}\label{avalue}
                a=\frac{4}{(3s_2-s_1)}
            \end{equation}

            There are three possibilities for $3s_2-s_1$.
            If $3s_2-s_1=1 \text{ or } 2$, then $a=4 \text{ or } 2$. But when $k$ is even, $a$ should be odd which is a contradiction. If $3s_2-s_1=4$, then $a=1$ is odd and contradicts \eqref{gcd}.\\
            By theorem\eqref{theorem1} and the necessary part of lemma\eqref{-1eigenlemma}, we can easily conclude that the multiplicity of the eigenvalue -1 is equal to $1$.
        \end{proof}
    \begin{lem}\label{1eigenlemma}
        For the graph $And(k)$, 1 is an eigenvalue if and only if $k\equiv3~mod~4$.
    \end{lem}
    \begin{proof}
    We can easily check that if $k=3~mod~4$, then $1$ is an eigenvalue of $And(k)$ by substituting $l=\frac{n}{4}$ in equation\eqref{eigenwhenkisodd} \\
    From the second part of the proof of lemma\eqref{-1eigenlemma}, it is easy to say that if 1 is an eigenvalue of $And(k)$, then $gcd(k+1,n)>1$ which implies that $k$ should be an odd integer.  Now suppose that $k\equiv 1~mod~4$, then $gcd(n,k+1)=2$ which gives $l=\frac{n}{2}$, which makes a contradiction that $x_{\frac{n}{2}}=-1$. Hence, $k\equiv 3~mod~4$.
    By theorem\eqref{theorem1} and the necessary part of lemma\eqref{1eigenlemma}, we can easily conclude that the multiplicity of the eigenvalue $1$ is equal to $1$.
    \end{proof}
    \begin{ex}
    For k=5, the graph And(5) is given in figure \ref{And(5)} and its adjacency matrix $A_{And(5)}$ is 
	\[A_{And(5)}=
	\begin{blockarray}{ccccccccccccccc}
		&0&1&2&3&4&5&6&7&8&9&10&11&12&13 \\
		\begin{block}{c(cccccccccccccc)}
			0&0&1&0&0&1&0&0&1&0&0&1&0&0&1\\
			1&1&0&1&0&0&1&0&0&1&0&0&1&0&0\\
			2&0&1&0&1&0&0&1&0&0&1&0&0&1&0\\
			3&0&0&1&0&1&0&0&1&0&0&1&0&0&1\\
			4&1&0&0&1&0&1&0&0&1&0&0&1&0&0\\
			5&0&1&0&0&1&0&1&0&0&1&0&0&1&0\\
			6&0&0&1&0&0&1&0&1&0&0&1&0&0&1\\
			7&1&0&0&1&0&0&1&0&1&0&0&1&0&0\\
			8&0&1&0&0&1&0&0&1&0&1&0&0&1&0\\
			9&0&0&1&0&0&1&0&0&1&0&1&0&0&1\\
			10&1&0&0&1&0&0&1&0&0&1&0&1&0&0\\
                11&0&1&0&0&1&0&0&1&0&0&1&0&1&0\\
                12&0&0&1&0&0&1&0&0&1&0&0&1&0&1\\
                13&1&0&0&1&0&0&1&0&0&1&0&0&1&0\\
		\end{block}
	\end{blockarray}
	\]
	The eigenvalues of the above matrix are
	\begin{align*}
		x_{0}&=5
		\\x_{1}&=x_{13}=0.356896
		\\x_{2}&=x_{12}=0.445042
		\\x_{3}&=x_{11}=0.692022
		\\x_{4}&=x_{10}=1.801938
		\\x_{5}&=x_{9}=-4.048917
		\\x_{6}&=x_{8}=-1.24698
		\\x_{7}&=-1
	\end{align*}
 In the above example, $And(5)$ has $k+\lceil\frac{k}{2}\rceil$=$5+\lceil\frac{5}{2}\rceil$=8 distinct eigenvalues.  Also the eigenvalues, $\lambda_{0}=x_{0}$, $\lambda_{1}=x_{k-1}=x_{2k}$, and $\lambda_{7}=x_{k}=x_{2k-1}=$ are 5, 1.801938 and -4.048917 respectively. Also, $And(5)$ has the eigenvalue $-1$ because $k$ is odd. However, in the above example $1$ is not an eigenvalue because  $k \not\equiv 3mod4$.
 \end{ex}
 \begin{ex}
 For k=3, the graph $And(3)$ is given in figure \ref{And(3)} and its adjacency matrix $A_{And(3)}$ is  
	\[A_{And(3)}=
	\begin{blockarray}{ccccccccc}
		&0&1&2&3&4&5&6&7 \\
		\begin{block}{c(cccccccc)}
			0&0&1&0&0&1&0&0&1\\
			1&1&0&1&0&0&1&0&0\\
			2&0&1&0&1&0&0&1&0\\
			3&0&0&1&0&1&0&0&1\\
			4&1&0&0&1&0&1&0&0\\
			5&0&1&0&0&1&0&1&0\\
			6&0&0&1&0&0&1&0&1\\
			7&1&0&0&1&0&0&1&0\\
		\end{block}
	\end{blockarray}
	\]
	The eigenvalues of the above matrix are
	\begin{align*}
		x_{0}&=3
		\\x_{1}&=x_{7}=0.4142
		\\x_{2}&=x_{6}=1
		\\x_{3}&=x_{5}=-2.4142
		\\x_{4}&=-1
	\end{align*}
 By examining the eigenvalues of the adjacency matrix for $And(3)$, we find that it has $3+\lceil\frac{3}{2}\rceil$=5 distinct eigenvalues.  Also, the eigenvalues  $\lambda_{0}=x_{0}$, $\lambda_{1}=x_{k-1}=x_{2k}$, and $\lambda_{7}=x_{k}=x_{2k-1}$ are $3$, $1$ and -$2.142$ respectively. Additionally, we observe that $-1$ and $1$ are eigenvalues due to the odd value of $k$, which follows the form $k\equiv 3mod4$.
 \end{ex}
 \begin{ex}
 For k=4, the graph, And(4) is given in figure \ref{And(4)} and its adjacency matrix $A_{And(4)}$ is 
	\[A_{And(4)}=
	\begin{blockarray}{cccccccccccc}
		&0&1&2&3&4&5&6&7&8&9&10 \\
		\begin{block}{c(ccccccccccc)}
			0&0&1&0&0&1&0&0&1&0&0&1\\
			1&1&0&1&0&0&1&0&0&1&0&0\\
			2&0&1&0&1&0&0&1&0&0&1&0\\
			3&0&0&1&0&1&0&0&1&0&0&1\\
			4&1&0&0&1&0&1&0&0&1&0&0\\
			5&0&1&0&0&1&0&1&0&0&1&0\\
			6&0&0&1&0&0&1&0&1&0&0&1\\
			7&1&0&0&1&0&0&1&0&1&0&0\\
			8&0&1&0&0&1&0&0&1&0&1&0\\
			9&0&0&1&0&0&1&0&0&1&0&1\\
			10&1&0&0&1&0&0&1&0&0&1&0\\
		\end{block}
	\end{blockarray}
	\]
	The eigenvalues of the above matrix are
	\begin{align*}
		x_{0}&=4
		\\x_{1}&=x_{10}=0.3728
		\\x_{2}&=x_{9}=0.5462
		\\x_{3}&=x_{8}=1.3979
		\\x_{4}&=x_{7}=-3.2287
		\\x_{5}&=x_{6}=-1.0882
	\end{align*}
By examining the eigenvalues of the adjacency matrix for $And(4)$, we find that it has $4+\lceil\frac{4}{2}\rceil$=6 distinct eigenvalues.  Also, the eigenvalues  $\lambda_{0}=x_{0}$, $\lambda_{1}=x_{k-1}=x_{2k}$, and $\lambda_{7}=x_{k}=x_{2k-1}$ are $4$, $1.3979$ and $-3.2287$ respectively. Since $k$ is even, $And(4)$ does not possess the eigenvalues $1$ and $-1$.
	
 \end{ex}
	
	\section{Conclusion}
	Our research has provided a comprehensive analysis of Andr\'asfai graphs, focusing on several key spectral properties. We have determined the second-largest eigenvalue of the Andr\'asfai graphs. This eigenvalue plays a crucial role in understanding the connectivity and structural properties of the graphs, shedding light on their spectral behavior.
Also, we have identified the smallest eigenvalue of the Andr\'asfai graph. This smallest eigenvalue holds important implications for various graph-theoretic properties, including connectivity, diameter, and expansion. Moreover, our research has elucidated the number of distinct eigenvalues exhibited by the Andr\'asfai graphs. 
    Furthermore, we have investigated the multiplicity of eigenvalues of  1 and -1 in the spectrum of the Andr\'asfai graphs. The determination of these multiplicities enhances our understanding of the symmetry and spectral symmetry of the graphs.
	\section*{Statements and Declarations}
	The authors declare that no funds, grants, or other support were received during the preparation of this manuscript.
	\section*{Availability of Data and Materials}
	Data sharing is not applicable to this article as no datasets were generated or analyzed during the current study.






\bigskip
\bigskip

{\footnotesize \pn{\bf Bharani Dharan K}\; \\ {Division of
Mathematics}, {School of Advanced Sciences}, {Vellore Institute of Technology}, {Chennai, India 600 127}\\
{\tt Email: bharanidharan.k2022@vitstudent.ac.in}\\

{\footnotesize \pn{\bf Radha S}\; \\ {Division of
Mathematics}, {School of Advanced Sciences}, {Vellore Institute of Technology}, {Chennai, India 600 127}\\
{\tt Email: radha.s@vit.ac.in}\\
\end{document}